\documentclass[11pt]{amsart}

\newtheorem{theorem}{Theorem}
\newtheorem{proposition}{Proposition}
\newtheorem{corollary}{Corollary}
\newtheorem{remark}{Remark}

\theoremstyle{definition}
\newtheorem{ex}{Example}

\DeclareMathOperator{\trace}{trace}
\DeclareMathOperator{\grad}{grad} \DeclareMathOperator{\rank}{rank}
 \DeclareMathOperator{\Div}{div}

\numberwithin{equation}{section}

\DeclareMathOperator{\ricci}{Ricci}
\DeclareMathOperator{\Vol}{Vol}

\begin{document}

\title[Biharmonic stress-energy tensor]{The biharmonic stress-energy tensor and the Gauss map}

\author{E.~Loubeau}
\author{S.~Montaldo}
\author{C.~Oniciuc}

\address{D{\'e}partement de Math{\'e}matiques \\
Laboratoire CNRS UMR 6205 \\
Universit{\'e} de Bretagne Occidentale \\
6, avenue Victor Le Gorgeu \\
CS 93837, 29238 Brest Cedex 3, France} \email{loubeau@univ-brest.fr}

\address{Universit\`a degli Studi di Cagliari\\
Dipartimento di Matematica e Informatica\\
Via Ospedale 72\\
09124 Cagliari, Italia} \email{montaldo@unica.it}

\address{Faculty of Mathematics\\ ``Al.I. Cuza'' University of Iasi\\
Bd. Carol I no. 11 \\
700506 Iasi, Romania} \email{oniciucc@uaic.ro}



\begin{abstract}
We consider the energy and bienergy functionals as variational
problems on the set of Riemannian metrics and present a study of the
biharmonic stress-energy tensor. This approach is then applied to
characterise weak conformality of the Gauss map of a submanifold.
Finally, working at the level of functionals, we recover a result of
Weiner linking Willmore surfaces and pseudo-umbilicity.
\end{abstract}

\maketitle

\section{Introduction}

The guiding principle of variational theory is that geometric
objects can be selected according to whether or not they minimize
certain functionals and, since Morse theory, critical points can
prove sufficient. Once this criterion chosen, the adequate
Euler-Lagrange equation will characterise maps particularly well
adapted to our geometric framework. However, roles can be reversed
and metrics can be viewed as variables and required to fit in with a
map and complete the picture. Other than the duality of these
approaches, the theory of general relativity has put metrics firmly
in centre stage and the characterisation of Einstein metrics as
(constrained) critical points of the total curvature has created a
new viewpoint on the usual functionals, in particular the various
energies defined for maps between manifolds.

Let $\phi:(M^m,g)\to (N^n,h)$ be a smooth map between Riemannian
manifolds, assume $M$ compact and define the {\it energy} of $\phi$
to be
$$
E(\phi)=\int_M e(\phi) \ v_g ,
$$
where $e(\phi) =\tfrac{1}{2} \vert d\phi\vert^2$ is (half) the
Hilbert-Schmidt norm.

\noindent Call a map {\it harmonic} if it is a critical point of
$E$, i.e. $\frac{d}{dt}\big{\vert}_{t=0}E(\phi_t)=0$, for any smooth
deformation $\{\phi_t\}$ of $\phi$. The corresponding Euler-Lagrange
equation characterizes harmonicity
\begin{eqnarray*}
\tau(\phi)&=g^{ij}\big(\frac{\partial^2 \phi^{\alpha}}{\partial
x^i\partial x^j}-^M\Gamma^k_{ij}\phi^{\alpha}_k
+^N\Gamma^{\alpha}_{\beta\sigma}\phi^{\beta}_i\phi^{\sigma}_j\big)\frac{\partial}{\partial
y^{\alpha}} =0,
\end{eqnarray*}
where $^M\Gamma^k_{ij}$ and $^N\Gamma^{\alpha}_{\beta\sigma}$ are
the Christoffel symbols of $g$ and $h$.\\
On non-compact manifolds, this equation serves as definition.

If $M$ is compact, the set $\mathcal{G}$ of Riemannian metrics on
$M$ is an infinite dimensional manifold and its tangent space at $g$
is identified with symmetric $(0,2)$-tensors:
$$
T_g\mathcal{G}=C(\odot^2T^{\ast}M).
$$
For a deformation $\{g_t\}$ of $g$ we denote
$\omega=\frac{d}{dt}\big{\vert}_{t=0}g_t\in T_g\mathcal{G}$.

Now, fix $\phi:M\to (N,h)$ and define the functional
$\mathcal{F}:\mathcal{G}\to {\mathbb R}$ by
$$
\mathcal{F}(g)=E(\phi),
$$
where $E(\phi)$ is computed with respect to the metrics $g$ and $h$.

Sanini obtained the Euler-Lagrange equation for $\mathcal{F}$.

\begin{theorem}[\cite{AS}]
Let $\phi:M\to (N,h)$ and assume that $M$ is compact, then
$$
\frac{d}{dt}\big{\vert}_{t=0}\mathcal{F}(g_t)=\frac{1}{2}\int_M
\langle \omega, e(\phi)g-\phi^{\ast}h \rangle \ v_g,
$$
so $g$ is a critical point of $\mathcal{F}$ if and only if the
stress-energy tensor $S=e(\phi)g-\phi^{\ast}h$ vanishes.
\end{theorem}

This naturally extends into a definition on non-compact domains and
Baird and Eells proved:

\begin{theorem}[\cite{PBJE}]
Let $\phi:(M,g)\to (N,h)$ be a map between Riemannian manifolds,
then:
$$
\Div S(X)=-\langle\tau(\phi),d\phi(X)\rangle, \quad \forall X\in
C(TM).
$$
Therefore, if $\phi$ is harmonic then $\Div S=0$.
\end{theorem}

The vanishing of $S$ is a strong condition:

\begin{theorem}[\cite{PBJE,AS}]
Let $\phi:(M,g)\to (N,h)$. Then $S=0$ if and only if either $m=2$
and $\phi$ is conformal, or $m>2$ and $\phi$ is constant.
\end{theorem}

Note that a homothetic transformation of the domain can render
$\mathcal{F}$ arbitrarily large or small, since
$\mathcal{F}(tg)=t^{\frac{m-2}{2}}\mathcal{F}(g)$, for a positive
constant $t$. To avoid this, impose $\Vol (M,g_t)=\Vol (M,g)$, i.e.
$\{g_t\}$ is an isovolumetric deformation, in this case $\omega$ is
orthogonal to $g$ as vectors in $T_{g}\mathcal{G}$, i.e.
$$
(\omega,g)=\int_M\langle \omega,g\rangle \ v_g=0,
$$
and $g$ is a critical point of $\mathcal{F}$ with respect to
isovolumetric deformations of $g$ if and only if $S=\lambda g$,
where $\lambda$ is a real constant.

\begin{theorem}[\cite{AS}]
Let $\phi:(M,g)\to (N,h)$. Then $S=\lambda g$ if and only if either
$m=2$ and $\phi$ is conformal, or $m>2$ and $\phi$ is a homothety.
\end{theorem}

\section{The biharmonic case}

Let $\phi:(M^m,g)\to (N^n,h)$ be a smooth map between Riemannian
manifolds, assume $M$ compact and define the {\it bienergy} of
$\phi$ by:
$$
E_2(\phi)=\frac{1}{2}\int_M \vert\tau(\phi)\vert^2 \ v_g.
$$

\noindent A map is called {\it biharmonic} if critical point of
$E_2$ and Jiang derived its Euler-Lagrange equation.

\begin{theorem}[\cite{GYJ1}]
Let $\phi:(M,g)\to (N,h)$ and assume $M$ compact. Then $\phi$ is
biharmonic if and only if
\begin{eqnarray*}
\tau_2(\phi)&=&-\Delta\tau(\phi)-\trace
R^N(d\phi\cdot,\tau(\phi))d\phi\cdot =0.
\end{eqnarray*}
\end{theorem}

\noindent In this paper we use the sign conventions
$\Delta\sigma=-\trace \nabla d\sigma$, $\sigma\in C(\phi^{-1}TN)$,
and $R(X,Y)Z=\nabla_X\nabla_YZ-\nabla_Y\nabla_XZ-\nabla_{[X,Y]}Z$.

Obviously, any harmonic map is biharmonic, therefore we are
interested in non-harmonic biharmonic maps, which we call {\it
proper biharmonic}.

Two examples of proper biharmonic maps are:
\begin{enumerate}
\item The inclusion ${\bf i}:{\mathbb S}^n(\frac{1}{\sqrt{2}})\to {\mathbb S}^{n+1}$ is proper
biharmonic
\item Let $\psi:M\to {\mathbb S}^n(\frac{1}{\sqrt{2}})$ be a harmonic map with
$e(\psi)$ constant. Then the composition map $\phi={\bf i}\circ\psi$
is proper biharmonic.
\end{enumerate}
For an account of biharmonic maps see~\cite{SMCO} and {\it The
bibliography of biharmonic maps}~\cite{BMBib}.

To a map $\phi:(M,g)\to (N,h)$, Jiang associates in~\cite{GYJ2} the
symmetric $(0,2)$ tensor:
\begin{eqnarray*}
S_2(X,Y)&=&\big(\tfrac{1}{2}\vert\tau(\phi)\vert^2+\langle
d\phi,\nabla\tau(\phi)\rangle\big)\langle X,Y\rangle -\langle
d\phi(X),\nabla_Y\tau(\phi)\rangle-\langle
d\phi(Y),\nabla_X\tau(\phi)\rangle
\end{eqnarray*}
and proved
\begin{equation}
\label{eq:div-S-2} \Div S_2(X)=-\langle\tau_2(\phi),d\phi(X)\rangle.
\end{equation}
Therefore, if $\tau_2(\phi)=0$ then $\Div S_2=0$.\\
As for harmonic maps, the expression of $S_2$ can be deduced from a
variational problem.

\begin{theorem}~\cite{ELSMCO}
Fix $\phi:M\to (N,h)$, assume $M$ compact and define
$\mathcal{F}_2:\mathcal{G}\to {\mathbb R}$ to be
$$
\mathcal{F}_2(g)=E_2(\phi),
$$
then
$$
\frac{d}{dt}\big\vert_{t=0}\mathcal{F}_2(g_t)=-\frac{1}{2}\int_M
\langle \omega, S_2\rangle \ v_g .
$$
So $g$ is a critical point of $\mathcal{F}_2$ if and only if
$S_2=0$.
\end{theorem}

\noindent From~\eqref{eq:div-S-2} we obtain:

\begin{proposition} If $\phi:(M,g)\to (N,h)$ is:
\begin{itemize}
\item [a)] a Riemannian immersion then $\Div S_2=0$ if
and only if $\tau_2(\phi)$ is normal.
\item [b)] a submersion (not necessarily Riemannian) then $\Div
S_2(\phi)=0$ if and only if $\tau_2(\phi)=0$.
\end{itemize}
\end{proposition}

\noindent This allows us to obtain new examples of proper biharmonic
maps.

\begin{proposition}~\cite{ELSMCO}
Let $\phi:(M,g)\to (N,h)$ be a submersion with basic tension field,
i.e. $\tau(\phi)=\xi\circ\phi$, $\xi\in C(TN)$, and $\xi$ Killing.
If $M$ is compact then $\phi$ is harmonic, while if $M$ is
non-compact then $\phi$ is proper biharmonic if and only if the norm
of $\xi$ is constant (non-zero).
\end{proposition}

\begin{ex} Let $(M^m,g)$ and $(N^n,h)$ be Riemannian manifolds and
$f\in C^{\infty}(M)$ a positive function. Consider the warped
product manifold $M\times_{f^2}N$, then the projection $\pi$ onto
the first term is a Riemannian submersion and $\tau(\pi)=n\grad (\ln
f)\circ\pi$. If $\ln f$ is an affine function on $M$ then $\grad
(\ln f)$ is a Killing vector field of constant norm and $\pi$ is
biharmonic.
\end{ex}

\begin{ex} For any vector field $\xi$, the tangent bundle
$TM$ can be endowed with a Sasaki-type metric such that the
canonical projection is a Riemannian submersion and
$\tau(\pi)=-(m+1)\xi\circ\pi$ (\cite{CO}). If $\xi$ is Killing of
constant norm then $\pi$ is biharmonic.
\end{ex}

If $\tau(\phi)=0$ then $S_2=0$ but the converse, i.e. $S_2=0$ (a
critical point of $\mathcal{F}_2$) implies $\tau(\phi)=0$ (an
absolute minimum of $\mathcal{F}_2$) is less straight-forward. Note
that, in general, $S_2=0$ does not imply harmonicity; for example,
the non-geodesic curve $\gamma(t)=t^3a$, $a\in {\mathbb R}^n$, has
$S_2=0$. Remember also that for harmonicity, when $m>2$, $S=0$
implies $\phi$ constant.

The vanishing of $S_2$ implies harmonicity in some situations
(confer~\cite{ELSMCO}):

\begin{enumerate}
\item curves parametrized by arc-length,
\item $\phi:(M^2,g)\to (N,h)$,
\item $\phi:(M^m,g)\to (N,h)$, $m>2$, and $\rank \phi\leq m-1$,
\item $\phi:(M^m,g)\to (N,h)$, $m>2$, and $\phi$ submersion,
\item $\phi:(M^m,g)\to (N,h)$, $m\neq 4$, $M$ compact (\cite{GYJ1}),
\item $\phi:(M^m,g)\to (N,h)$ Riemannian immersion, $m\neq 4$.
\end{enumerate}

Dimension $4$ plays a special role for the domain manifold, as we
can see from the followings

\begin{theorem}[\cite{GYJ2}]
Let $\phi:(M^4,g)\to (N,h)$ be a non-minimal Riemannian immersion,
then $S_2=0$ if and only if $\phi$ is pseudo-umbilical.
\end{theorem}

\noindent To generalize this result, we have to consider conformal
immersions:

\begin{proposition}~\cite{ELSMCO}
Let $\phi:(M^4,g=e^{2\rho}\phi^{\ast}h)\to (N,h)$ be a conformal
immersion, $M$ compact. Then $S_2=0$ if and only if $\rho$ is
constant and $\overline{\phi}:(M^4 ,\phi^{\ast}h)\to (N,h)$ is
pseudo-umbilical.
\end{proposition}

\begin{proposition}~\cite{ELSMCO}
Let $\phi:(M^4,g)\to (N^4,h)$ be a local diffeomorphism, i.e. $\rank
\phi=4$, $M$ compact. Then $S_2=0$ if and only if $\tau(\phi)=0$.
\end{proposition}

\begin{proposition}~\cite{ELSMCO}
Let $\phi:(M^4,g)\to (N,h)$ be a map such that $\rank \phi\leq 3$.
Then $S_2=0$ if and only if $\tau(\phi)=0$.
\end{proposition}

Then we consider isovolumetric deformations of the domain metric:

\begin{theorem}~\cite{ELSMCO}
Let $\phi:(M^m,g)\to (N,h)$ be a Riemannian immersion. Then
$S_2=\lambda g$ if and only if either $m=4$ and $\phi$ is
pseudo-umbilical, or $m\neq 4$ and $\phi$ is pseudo-umbilical with
$\vert\tau(\phi)\vert$ constant.
\end{theorem}

We end this section with the study of the behaviour of $S_2$ under
conformal changes of the domain metric.

\begin{proposition}
Consider $\phi:(M^m,g)\to (N^n,h)$,
$\tilde{\phi}:(M,\tilde{g}=tg)\to (N,h)$,
$\phi=\tilde{\phi}\circ{\bf 1}$, where ${\bf 1}:(M,g)\to
(M,\tilde{g})$ is the identity map and $t$ is a positive constant.
Then $ \tilde{S_2}=\frac{1}{t}S_2, $ therefore $\tilde{S_2}=0$ if
and only if $S_2=0$.
\end{proposition}

\noindent For surfaces we get:

\begin{proposition}
Let $\phi:(M^2,g)\to (N^n,h)$ and
$\tilde{\phi}:(M,\tilde{g}=e^{2\rho}g)\to (N,h)$,
$\phi=\tilde{\phi}\circ{\bf 1}$, $\rho\in C^{\infty}(M)$:
\begin{itemize}
\item[a)] $\tilde{S_2}=0$ if and only if $S_2=0$ and, in this case,
the maps are harmonic.
\item[b)] if $\langle\tau(\phi),d\phi(X)\rangle=0$, $\forall X\in C(TM)$,
then $ \tilde{S_2}=e^{-2\rho}S_2 $.
\end{itemize}
\end{proposition}

\noindent For domains of higher dimension we obtain two ``rigidity''
results:

\begin{proposition}
Let $M^m$ be compact, $m>2$, $m\neq 4$. Consider
$\phi:(M^m,g)\to(N^n,h)$ such that
$\langle\tau(\phi),d\phi(X)\rangle=0$, $\forall X\in C(TM)$ and
$\tilde{\phi}:(M,\tilde{g}=e^{2\rho}g)\to (N,h)$. Then
$\tilde{S_2}=0$ if and only if $d\phi(\grad\rho) =0$ and $S_2=0$,
and both maps must then be harmonic. \\
When $\phi$ is a Riemannian immersion, $\tilde{S_2}=0$ if and only
if $\rho$ is constant and $S_2=0$.
\end{proposition}

\begin{proposition}~\cite{ELSMCO}
Let $\phi:(M^4,g)\to (N^n,h)$ be a non-minimal Riemannian immersion
and assume that $M$ is compact. Let
$\tilde{\phi}:(M,\tilde{g}=e^{2\rho}g)\to (N,h)$, then
$\tilde{S_2}=0$ if and only if $\rho$ is constant and $S_2=0$. In
this case $\phi$ is pseudo-umbilical.
\end{proposition}

\section{The tensor $S_2$ and the Gauss map}

Let $M^m$ be an oriented submanifold of ${\mathbb R}^n$, $p\in M$ an
arbitrary point and $\{X_i\}_{i=1}^m$ a positive oriented geodesic
basis centered around $p$. On a neighbourhood $U$ of $p$, the Gauss
map associated to $M$ can be written:
\begin{align*}
G&:M\to G(n,m) \\G(q)&=X_1(q)\wedge\ldots\wedge X_m(q), \quad
\forall q\in U.
\end{align*}
Since
$$
dG_q(X_i)=\sum_{j=1}^m X_1(q)\wedge\ldots\wedge
X_{j-1}(q)\wedge\big(\nabla^0_{X_i}X_j\big)(q)\wedge
X_{j+1}(q)\wedge\ldots\wedge X_m(q),
$$
where $\nabla^0$ is the canonical connection on ${\mathbb R}^n$, at
$p$ we have:
$$
dG_p(X_i)=\sum_{j=1}^m X_1(p)\wedge\ldots\wedge X_{j-1}(p)\wedge
B_p(X_i,X_j)\wedge X_{j+1}(p)\wedge\ldots\wedge X_m(p),
$$
where $B$ denotes the second fundamental form of $M$.

\noindent Complete $\{X_i(p)\}_{i=1}^m$ into an orthonormal basis
$\{X_{\alpha}(p)\}_{\alpha=1}^n$ of ${\mathbb R}^n$. Let $\alpha\in
\{1,\ldots,n\}$ and $a\in\{m+1,\ldots,n\}$, then:
$$
B_p(X_i,X_j)=\sum_{a} b_{ij}^a(p)X_a(p),
$$
and
$$
dG_p(X_i)=\sum_a\sum_jb_{ij}^a(p)X_1(p)\wedge\ldots\wedge
X_{j-1}(p)\wedge X_a(p)\wedge X_{j+1}(p)\wedge\ldots\wedge X_m(p).
$$
Now, the $m-$subspace $X_1(p)\wedge\ldots\wedge X_{j-1}(p)\wedge
X_a(p)\wedge X_{j+1}(p)\wedge\ldots\wedge X_m(p)$, can be identified
with $X_j^{\ast}(p)\otimes X_a(p)$ (\cite{JELL}), so
$$
dG_p(X_i)=\sum_a\sum_jb_{ij}^a(p)X_j^{\ast}(p)\otimes X_a(p).
$$
\noindent The canonical metric $g_{can}$ on $G(n,m)$ is defined by
requiring that
$$
\{X_j^{\ast}(p)\otimes X_a(p) : j=1,\ldots,m, \ a=m+1,\ldots, n\}
$$
is an orthonormal basis of $T_{G(p)}G(n,m)$. By direct computation,
we obtain:
$$
g_{can}(dG_p(X_i),dG_p(X_k))=\sum_j \langle B_p(X_i,X_j),
B_p(X_k,X_j)\rangle,
$$
where $\langle,\rangle$ is the canonical metric on ${\mathbb R}^n$.
By the Gauss Lemma
\begin{eqnarray*}
g_{can}(dG_p(X_i),dG_p(X_k))=-\ricci_p(X_i,X_k)+m\langle
H(p),B_p(X_i,X_k)\rangle,
\end{eqnarray*}
where $H$ is the mean curvature vector field. Therefore
$$
(G^{\ast}g_{can})(p)=m\langle H(p),B_p\rangle-\ricci_p.
$$
Now
\begin{eqnarray*}
S^G&=&e(G)g-G^{\ast}g_{can}
=(\ricci-\frac{r}{2}g)+\frac{m^2}{2}\vert H\vert^2g-m\langle
H,B\rangle \\
&=&(\ricci-\frac{r}{2}g)+\frac{1}{2}\vert\tau({\bf
i})\vert^2g-\langle\tau({\bf i}),\nabla d{\bf i}\rangle \\
&=&(\ricci-\frac{r}{2}g)-\frac{1}{2}S_2^{\bf
i}+\frac{1}{4}\vert\tau({\bf i})\vert^2g,
\end{eqnarray*}
where $g=\langle,\rangle$, ${\bf i}$ is the canonical inclusion of
$M$ in ${\mathbb R}^n$ and $r=\trace \ricci$ is the scalar
curvature.

\begin{proposition}
Assume $M^2$ is an orientable surface in ${\mathbb R}^n$, then the
following conditions are equivalent:
\begin{itemize}
\item [a)] $S^G=0$,
\item [b)] $G$ is weakly conformal,
\item [c)] $M^2$ is pseudo-umbilical,
\item [d)] $S^{\mathbf i}_2=\frac{1}{2}\vert\tau({\bf i})\vert^2g$.
\end{itemize}
\end{proposition}

\begin{proposition}
Assume that $m>2$, then any two of the following statements implies
the third:
\begin{itemize}
\item[a)] $S^{\bf i}_2=fg$, where $f\in C^{\infty}(M)$,
\item [b)] M is Einstein,
\item [c)] $G$ is weakly conformal.
\end{itemize}
\end{proposition}

\begin{remark} We have:
\begin{itemize}
\item [a)]
if $S^{\bf i}_2=fg$ and $G$ is weakly conformal then $S^{\bf
i}_2=\frac{4-m}{2m}\vert\tau({\bf i})\vert^2g$,
$G^{\ast}g_{can}=\frac{2}{m}e(G)g$ and
$$
\ricci=\frac{\vert\tau({\bf i})\vert^2-2e(G)}{m}g,
$$
i.e. $M$ is Einstein. Moreover, in this case, $r=\vert\tau({\bf
i})\vert^2- 2e(G)$ must be constant.
\item[b)] if $S^{\bf i}_2=fg$ and $\ricci=cg$, $c$ constant, then $G$ is
weakly conformal and
$$
e(G)=\frac{\vert\tau({\bf i})\vert^2-mc}{2}.
$$
Moreover, in this case, $\vert\tau({\bf i})\vert^2-mc\geq 0$, and,
if $M$ has constant mean curvature then $G$ is homothetic. We
conclude that if $M^m$, $m>2$, is an Einstein pseudo-umbilical
submanifold of ${\mathbb R}^n$, with constant mean curvature when
$m\neq 4$, then its Gauss map is homothetic.
\end{itemize}
\end{remark}

Since $\Div(\ricci-\frac{r}{2}g)=0$ we re-obtain Jiang's result:

\begin{theorem}[\cite{GYJ2}]
Let $M^m$ be an oriented submanifold of ${\mathbb R}^n$. Then the
tensors $S^G$ and $S^{\mathbf i}_2$ are related by
$$
\Div S^G+\tfrac{1}{2}\Div S^{\bf i}_2-\tfrac{1}{4}d(\vert\tau({\bf
i})\vert^2)=0.
$$
\end{theorem}

Since, Ruh-Vilms proved in~\cite{R-V} that $G$ is harmonic if and
only if the mean curvature vector field is parallel, we conclude:

\begin{corollary}Let $M^m$ be an oriented submanifold of ${\mathbb
R}^n$, then:
\begin{itemize}
\item [a)] if $M$ has constant mean curvature, then $\Div S^{\bf i}_2=0$
if and only if $\Div S^G=0$,
\item[b)] if $G$ is harmonic then $\Div S^{\bf i}_2=0$.
\end{itemize}
\end{corollary}

\section{On a result of Weiner}

Inspired by the above technique on the Gauss map, we conclude with a
result on Willmore surfaces in ${\mathbb R}^n$ due to Weiner
in~\cite{JLW}.

Let $\phi:(M,g)\to {\mathbb R}^n$ be a Riemannian immersion, i.e.
$g=\phi^{\ast}\langle,\rangle$, assume $M$ oriented. We have
\begin{eqnarray*}
G^{\ast}g_{can}&=&m\langle H,B(\cdot,\cdot)\rangle-\ricci
=\langle\tau(\phi),\nabla d\phi(\cdot,\cdot)\rangle-\ricci,
\end{eqnarray*}
and
$$
e(G)=\frac{1}{2}m^2\vert H\vert^2-\frac{1}{2}r.
$$
Assume $m=2$, therefore
$$
e(G)=2\vert H\vert^2-K,
$$
where $K$ is the Gaussian curvature of $(M,g)$, and integrating,
$$
\int_M e(G) \ v_g=2\int_M \vert H\vert^2 \ v_g - 2\pi\chi(M).
$$

Consider a one-parameter family of immersions $\{\phi_t\}$,
$\phi_0=\phi$ such that $\phi_t:(M,g_t)\to {\mathbb R}^n$ is a
Riemannian immersion, i.e. $g_t=\phi_t^{\ast}\langle,\rangle$.  All
previous formulas hold for $\phi_t$, so, for any $t$:
$$
\int_M e(G_t) \ v_{g_t}=2\int_M\vert H_t\vert^2 \ v_{g_t} -
2\pi\chi(M).
$$
The right-hand side consists of the Willmore functional plus the
Euler-Poincar{\'e} characteristic, a topological invariant. Compute
$$
W=\frac{d}{dt}\Big\vert_{t=0}\int_M e(G_t) \ v_{g_t}
=2\frac{d}{dt}\Big\vert_{t=0}\int_M\vert H_t\vert^2 \ v_{g_t} .
$$
Put $h=g_{can}$, then:
\begin{eqnarray*}
2W&=&\frac{d}{dt}\Big\vert_{t=0}\int_M
g^{ij}(x,t)G^{\alpha}_i(x,t)h_{\alpha\beta}(G(x,t))G^{\beta}_j(x,t)
\ v_{g_t(x)}  \\
&=&\int_M \frac{\partial g^{ij}}{\partial
t}(x,0)G^{\alpha}_i(x)h_{\alpha\beta}(G(x))G^{\beta}_j(x) \ v_{g} \\
&&+\int_M
g^{ij}(x)\frac{d}{dt}\big\vert_{t=0}\big\{G^{\alpha}_i(x,t)h_{\alpha\beta}(G(x,t))G^{\beta}_j(x,t)\big\}
\ v_{g} \\
&&+\int_M
g^{ij}(x)G^{\alpha}_i(x)h_{\alpha\beta}(G(x))G^{\beta}_j(x) \
\frac{d}{dt}\big\vert_{t=0}v_{g_t(x)}.
\end{eqnarray*}
Let
$$
W_1= \int_M
g^{ij}(x)\frac{d}{dt}\big\vert_{t=0}\big\{G^{\alpha}_i(x,t)h_{\alpha\beta}(G(x,t))G^{\beta}_j(x,t)\big\}
\ v_{g}.
$$
so
\begin{eqnarray*}
2W&=&\int_M \frac{\partial g^{ij}}{\partial
t}(x,0)G^{\alpha}_i(x)h_{\alpha\beta}(G(x))G^{\beta}_j(x) \ v_{g} \\
&&+\int_M \big(4\vert H\vert^2-r\big) \
\frac{d}{dt}\big\vert_{t=0}v_{g_t(x)} \\
&&+ W_1.
\end{eqnarray*}

\noindent Recall that
$$
\frac{\partial g^{ij}}{\partial t}(x,0)=-g^{ik}g^{jl}\omega_{kl}
\quad \hbox{and} \quad \frac{d}{dt}\big\vert_{t=0}v_{g_t(x)}=\langle
\frac{1}{2}g,\omega\rangle v_g.
$$
Replacing we obtain:
\begin{eqnarray*}
2W&=&-\int_M\langle\omega, G^{\ast}h\rangle \ v_g + \int_M
\big(4\vert H\vert^2-r\big)\langle \frac{1}{2}g,\omega\rangle \ v_g
+W_1
\\
&=& \int_M\langle 2\vert H\vert^2 g-2\langle
H,B(,)\rangle,\omega\rangle \ v_g +W_1.
\end{eqnarray*}
Clearly if $G$ is harmonic (so $W_1=0$) and $M^2$ pseudo-umbilical
in $\mathbb{R}^n$ (i.e. $\vert H\vert^2 g-\langle H,B(,)\rangle=0$)
then it is Willmore.\\
To obtain a (partial) converse, we first establish the link between
$\omega$ and $V=\frac{d}{dt}\big\vert_{t=0}\phi_t$, which we assume
normal. Since $\omega = \frac{\partial g_{ij}}{\partial t}(x,0)
dx^{i} dx^{j}$ and $g_{ij}(x,t) = \sum_{\alpha=1}^{n}
\Phi^{\alpha}_{i}(x,t)\Phi^{\alpha}_{j}(x,t), \, i,j=1,2:$
\begin{align*}
\frac{\partial g_{ij}}{\partial t}(x,0) &=
2\sum_{\alpha}\frac{\partial^2 \Phi^{\alpha}}{\partial x^i
\partial t}(x,0) \phi^{\alpha}_{j}(x)
=2\sum_{\alpha}\frac{\partial V^{\alpha}}{\partial
x^i}(x)\phi^{\alpha}_{j}(x) \\
&=2\langle\nabla_{\partial_{i}} V, d\phi (\partial_{j})\rangle
=-2\langle V,\nabla_{\partial_{i}} d\phi (\partial_{j})\rangle\\
&=-2\langle V,\nabla d\phi (\partial_{i},\partial_{j}) -
d\phi(\nabla_{\partial_{i}}\partial_{j})\rangle =-2\langle
V,B(\partial_{i},\partial_{j})\rangle
\end{align*}
hence $\omega = -2  V.B$, where $(V.B)(X,Y)= \langle
V,B(X,Y)\rangle$. Therefore
\begin{align*}
\langle |H|^2 g - H.B , \omega\rangle &= -2 \langle
|H|^2 g - H.B ,  V.B \rangle \\
&= -2 |H|^2 \langle g, V.B\rangle + 2 \langle H.B , V.B \rangle
\end{align*}
but
$$\langle g, V.B\rangle = \sum_{i} \langle
V,B(X_{i},X_{i})\rangle = m \langle V,H\rangle = 2 \langle
V,H\rangle$$ and
\begin{align*}
\langle H.B , V.B \rangle &= \sum_{i,j}
\langle H,B(X_{i},X_{j})\rangle\langle V,B(X_{i},X_{j})\rangle \\
&=\sum_{i,j}
\Big(\sum_{a}  H^{a}B^{a}(X_{i},X_{j})\Big)\Big(\sum_{b} V^{b}B^{b}(X_{i},X_{j})\Big) \\
&= \sum_{b} \Big(\sum_{i,j,a}
H^{a}B^{a}(X_{i},X_{j})B^{b}(X_{i},X_{j})\Big)V^{b} ,
\end{align*}
where $a,b=3,\dots,n$. On the other hand, the contraction $\langle
H.B,B \rangle$ is the normal vector field defined by:
\begin{align*}
 \sum_{i,j} \langle
H,B(X_{i},X_{j}) \rangle B(X_{i},X_{j})  = \sum_{b}\sum_{i,j}\Big(
\sum_{a} H^{a}B^{a}(X_{i},X_{j})\Big) B^{b}(X_{i},X_{j}) \eta^{b} ,
\end{align*}
where $\{\eta^{b}\}$ is a normal frame, therefore:
$$
\langle\langle H.B,B \rangle ,V\rangle
=\sum_{b}\Big( \sum_{i,j} \Big(
\sum_{a}H^{a}B^{a}(X_{i},X_{j})\Big)B^{b}(X_{i},X_{j})\Big) V^{b}.
$$
Hence $\langle H.B , V.B \rangle =\langle\langle H.B,B \rangle
,V\rangle$ and
\begin{align*}
\langle |H|^2 g - H.B , \omega\rangle &= -4 |H|^2 \langle H,V
\rangle + 2\langle\langle H.B ,B \rangle
,V \rangle \\
&= \langle -4 |H|^2 H + 2\langle H.B ,B \rangle ,V \rangle.
\end{align*}
This shows that if we assume $G$ harmonic and $M^2$ Willmore then
$$
\int_M\langle -4 |H|^2 H + 2\langle H.B ,B \rangle ,V \rangle \
v_g=0
$$
for all normal variations $V$, as required by the Willmore problem,
hence $-4 |H|^2 H + 2\langle H.B ,B \rangle=0$. To conclude we need
to show that $-2 |H|^2 H + \langle H.B ,B \rangle =0$, or , since
$H= \tfrac{1}{2}\langle g,B\rangle$, $\langle - |H|^2 g +  H.B ,B
\rangle =0$, implies $- |H|^2 g + H.B =0$, i.e. $M^2$ is
pseudo-umbilical.\\
Decompose $B$ into its trace and trace-less parts: $B = H\otimes g +
S$, with $\trace S =0$, then $M^2$ is pseudo-umbilical if and only
if $ S .H  =0$ (umbilical being $S=0$). Then:
\begin{align*}
0&=\langle - |H|^2 g +  H.B ,B \rangle  \\
&= \langle - |H|^2 g +  H. (H\otimes g + S) ,H\otimes g
+ S \rangle \\
&=\langle - |H|^2 g + |H|^2 g +  H.S ,H\otimes g
+ S \rangle \\
&=\langle H,\trace S\rangle H + \sum_{i,j} \langle H,
S(X_{i},X_{j})\rangle S(X_{i},X_{j})
\end{align*}
therefore $\sum_{i,j} \langle H, S(X_{i},X_{j})\rangle
S(X_{i},X_{j}) =0$ and taking its inner-product with $H$, yields $
S.H  =0$.

Therefore we recover (part of) Weiner's result:

\begin{theorem}[\cite{JLW}]
Let $\phi:M^2\to {\mathbb R}^n$ be a Riemannian immersion of a
compact oriented surface into ${\mathbb R}^n$, such that its Gauss
map is harmonic. Then $M^2$ is a Willmore surface if and only if it
is pseudo-umbilical.
\end{theorem}

\begin{remark}
Recall Chen and Yano's result~\cite{BYC}: A submanifold of
$\mathbb{R}^n$ is pseudo-umbilical with parallel mean curvature
vector field if and only if it is minimal in a hypersphere of
${\mathbb R}^n$. So a minimal surface of $\mathbb{S}^{n-1}$ is a
Willmore surface of $\mathbb{R}^n$.
\end{remark}

\begin{remark}
The only compact oriented Riemannian immersed Willmore surface in
$\mathbb{R}^3$ of constant mean curvature is the sphere.
\end{remark}

\noindent {\it Acknowledgements}: The third author was partially
supported by the Grant At, 191/2006, C.N.C.S.I.S., Romania.


\begin{thebibliography}{99}

\bibitem{PBJE} Baird P. and Eells J., {\it A conservation law for harmonic
maps}, In: Geometry Symposium, Utrecht 1980, {\it Lecture Notes in
Math.}, 894, Springer 1981, pp 1--25.

\bibitem{BYC} Chen B.Y., {\it Geometry of Submanifolds},
Marcel Dekker, Inc., New York, 1973.

\bibitem{JELL} Eells J. and Lemaire L., {\it A report on harmonic
maps}, Bull. London Math. Soc. {\bf 10} (1978) 1--68.

\bibitem{GYJ1} Jiang G.Y., {\it 2-harmonic maps and their first and second
variational formulas}, Chinese Ann. Math. Ser. A {\bf 7} (1986)
389--402.

\bibitem{GYJ2} Jiang G.Y., {\it The conservation law for 2-harmonic
maps between Riemannian manifolds}, Acta Math. Sinica {\bf 30}
(1987) 220--225.

\bibitem{BMBib} The Bibliography of Biharmonic Maps,
{\tt http://beltrami.sc.unica.it/biharmonic/}

\bibitem{ELSMCO} Loubeau E., Montaldo S. and Oniciuc C., {\it The
stress-energy tensor for biharmonic maps}, {\tt
arXiv:math.DG/0602021}.

\bibitem{SMCO} Montaldo S. and Oniciuc C., {\it A short survey on
biharmonic maps between Riemannian manifolds}, {\tt
arXiv:math.DG/0510636}.

\bibitem{CO} Oniciuc C., {\it Biharmonic maps between Riemannian
manifolds}, An. Stiint. Univ. Al.I.~Cuza Iasi Mat. (N.S.) {\bf 48}
(2002) 237--248.

\bibitem{R-V} Ruh E. and Vilms J., {\it The tension field of the Gauss
map}, Trans. Amer. Math. Soc. {\bf 149} (1970) 569--573.

\bibitem{AS} Sanini A., {\it Applicazioni tra variet{\`a} riemanniane con
energia critica rispetto a deformazioni di metriche}, Rend. Mat.
{\bf 3} (1983) 53--63.

\bibitem{JLW} Weiner J., {\it On a problem of Chen, Willmore, et al}, Indiana Univ. Math.
J. {\bf 27} (1978) 19--35.

\end{thebibliography}
\end{document}